\documentclass[a4paper,11pt]{article}
\usepackage{latexsym}
\usepackage{amssymb}
\usepackage{amsfonts}
\usepackage{amsmath}
\usepackage{indentfirst}
\usepackage{graphicx}
\usepackage{epsfig}
\newtheorem{prop}{Proposition}[section]
\newtheorem{teor}{Theorem}[section]

\newtheorem{cor}{Corollary}[section]

\newcommand{\cvd}{\quad $\blacksquare$\bigskip}
\date{}

\author{Antonio Bernini\thanks{Dipartimento di
Sistemi e Informatica, viale Morgagni 65, 50134 Firenze, Italy
{\tt {bernini,ferrari,pinzani}@dsi.unifi.it}}\and Luca Ferrari$^\dag$\and Renzo Pinzani$^\dag$}

\title{Enumeration of some classes of words avoiding two generalized patterns of length
three\footnote{This work has been partially supported by MIUR
project: \emph{Automi e linguaggi formali: aspetti matematici e
applicativi}.}}

\begin{document}

\maketitle

\begin{abstract}
The method we have applied in \cite{BFP} to count pattern avoiding
permutations is adapted to words. As an application, we enumerate
several classes of words simultaneously avoiding two generalized
patterns of length 3.
\end{abstract}

\section{Introduction}

In the present work we deal with pattern avoidance on words. This
topic has first appeared in \cite{R}, and has been systematically
developed by Burstein in \cite{B}. Subsequently several authors
have studied this kind of matters, and in particular in \cite{BM}
exact formulas and/or generating functions for the number of words
avoiding a single generalized pattern of length 3 have been found.
Here we use a general method to count words on a totally ordered
alphabet avoiding a set of generalized patterns of length 3 of
type $(1,2)$ (i.e., having a dash between the first and the second
element). Our approach consists of inserting a letter at the end
of a given word of length $n$, thus obtaining a word of length
$n+1$. We perform this operation in such a way that part of the
preceding letters may have to be renamed. The choice of the letter
to be inserted depends on the patterns to be excluded. Obviously,
the above mentioned insertion technique, if applied to a word
avoiding the requested patterns, produces words in which the only
occurrence of a forbidden pattern could involve the newly inserted
element. Moreover, the particular type of the patterns to be
excluded allows to easily control that words of increasing length
will be correctly generated. The first appearance of this
technique goes back to \cite{BFP}; later, in \cite{E} the author
used it to count several classes of generalized pattern avoiding
permutations.

\bigskip

For the sake of clearness, we recall the above mentioned technique
in the framework of permutations. Given a permutation $\pi\in S$
(where $S$ denotes the whole symmetric group), it can be
represented in a path-like form as in the following figure:
\begin{center}
\includegraphics[scale=0.35]{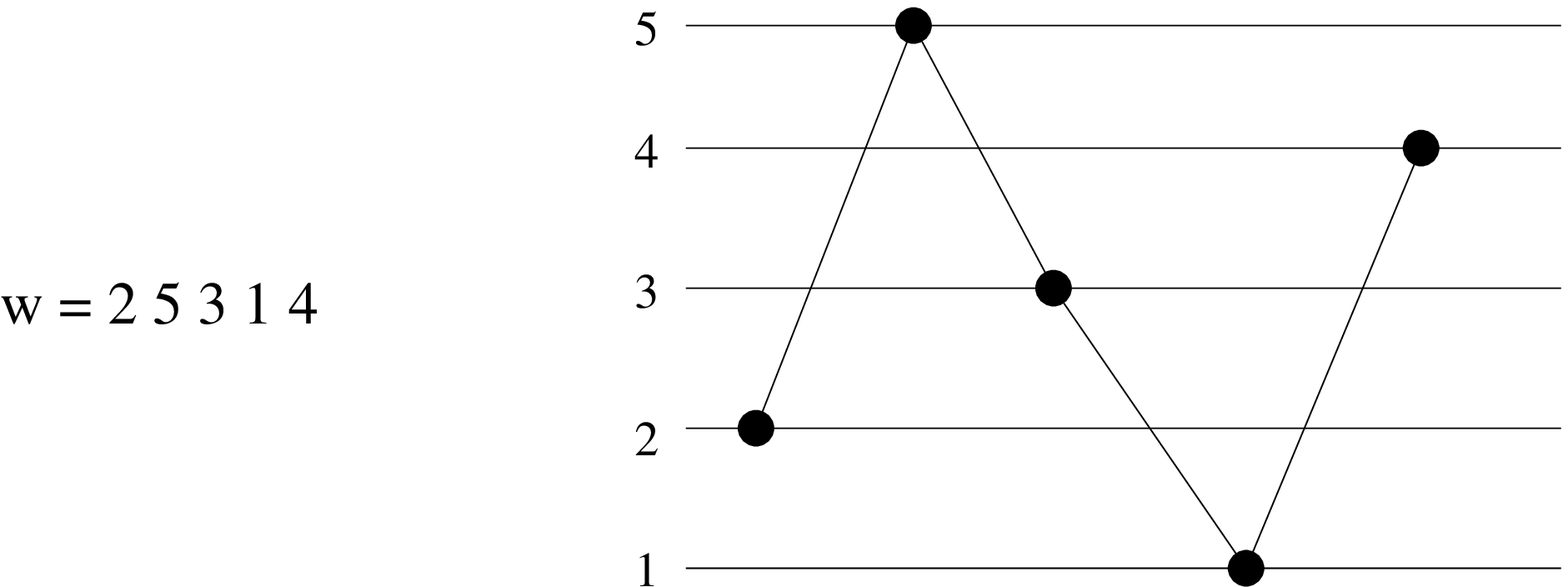}
\end{center}

Each entry of $\pi$ is a ``node" on the line corresponding to its
value. Every pair of adjacent nodes $\pi_i\ \pi_{i+1}$ is linked
by an ascending  or descending segment depending on whether
$\pi_i<\pi_{i+1}$ (i. e., the pair is an ``ascent") or
$\pi_i>\pi_{i+1}$ (i. e., the pair is a ``descent"). If $\pi \in
S_n$ (which is the set of all permutations of length $n$), the $n$
horizontal lines of the permutation divide the plane into $n+1$
regions, numbered $1$ to $n+1$ from bottom to top. Therefore, we
can obtain $n+1$ permutations belonging to $S_{n+1}$ starting from
$\pi$, by inserting a new node in each of these regions (see
figure below) and renaming the entries of the new permutation
$\pi' \in S_{n+1}$ according to the following \emph{renaming
rule}: if we insert the node into region $i$, then
\begin{enumerate}
    \item $\pi'_{n+1}=i$
    \item for $j=1,\ldots,n$
        \begin{enumerate}
            \item if $\pi_j<i$ then $\pi'_j=\pi_j$;
            \item otherwise $\pi'_j=\pi_j+1$.
        \end{enumerate}
\end{enumerate}
        \begin{center}
        \includegraphics[scale=0.35]{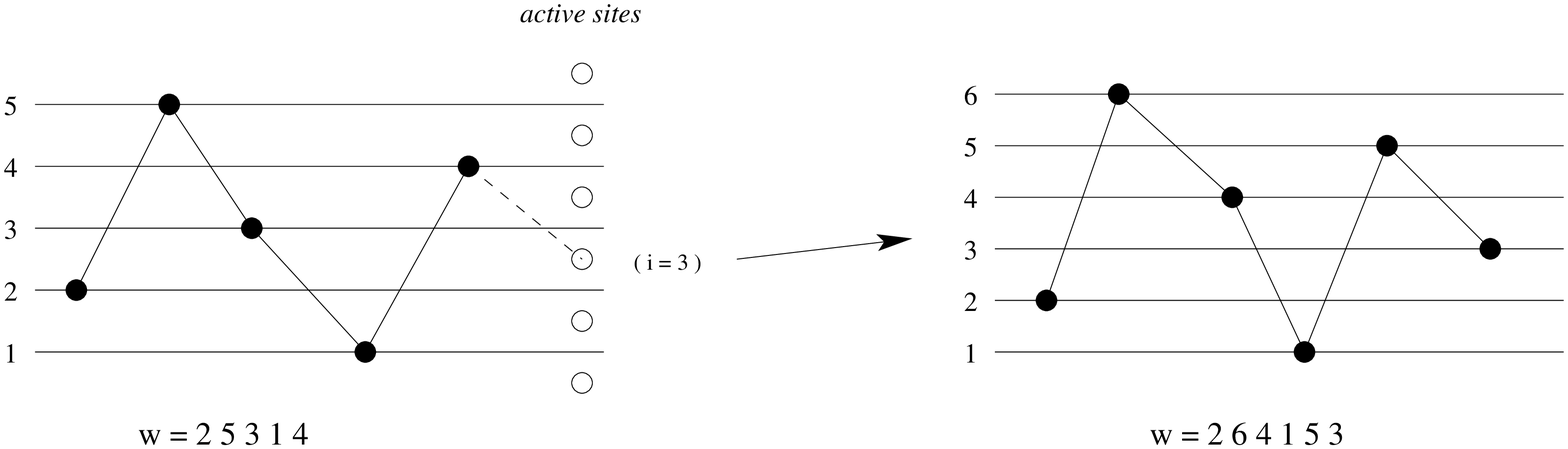}
        \end{center}

Now, we briefly recall the notion of (generalized) pattern
avoidance. If $\pi =\pi_1 \cdots \pi_n$ and $\sigma =\sigma_1
\cdots \sigma_k$ are two permutations of $S_n$ and $S_k$,
respectively, then $\pi$ \emph{avoids} the pattern $\sigma$ if
there are no indexes $i_1 <i_2 <\cdots <i_k$ such that
$\pi_{i_1}\pi_{i_2}\cdots \pi_{i_k}$ is in the same relative order
as $\sigma_{i_1}\sigma_{i_2}\cdots \sigma_{i_k}$. The subset of
$\sigma$-avoiding permutations of $S_n$ is denoted $S_n (\sigma)$.
A \emph{generalized} pattern $\tau$ \cite{BS} is a pattern with
some dashes inserted, so that two consecutive elements of $\tau$
are adjacent if there is no dash between them (e.g. $261-4-35$ is
a generalized pattern of length 6). A permutation $\pi$
\emph{contains} the generalized pattern $\tau$ if the elements of
$\pi$ corresponding to the elements of $\tau$ are in the same
relative order and any two elements of $\pi$ corresponding to two
adjacent elements of $\tau$ are adjacent in $\pi$ as well. A
permutation $\pi$ \emph{avoids} a generalized pattern $\tau$ if
$\pi$ does not contain $\tau$. For instance, the permutation
$41325$ contains $321$ but not $32-1$. In general, if $T$ is a set
of (generalized) patterns, $S_n (T)$ denotes the set of
permutations of $[n]=\{1,2,\ldots,n\}$ avoiding each pattern of
$T$.

\bigskip

If we want to generate all the permutations in $S_{n+1}(T)$
avoiding certain generalized patterns, then the regions where we
can insert the new node form a subset of all the $n+1$ possible
regions, whose elements are called \emph{active sites}. The
insertion of the new node to generate a new permutation leads to
an important consideration: if $\pi\in S_n(T)$, then $\pi'$ ($\in
S_{n+1}$), obtained from $\pi$ inserting the last node in some
region, does not contain the patterns specified in $T$ in its
entries $\pi'_j$ with $j=1,\ldots,n$, otherwise $\pi$ itself would
contain some pattern of $T$. Therefore, we can decide if a region
$i$ is an active site or not simply by checking those generalized
patterns the last node is involved in.

For the enumeration of a class of generalized pattern avoiding
permutations (according to their length), the following
general strategy can be considered:%
%
\begin{itemize}
    \item find the active sites for the permutations of $S_n(T)$
    among the $n+1$ possible regions;
    \item describe the generation of the permutations of
    $S_{n+1}(T)$ and encode the construction with a succession rule;
    \item obtain the generating function of the sequence
    enumerating $|S_n(T)|$ from the succession rule with analytic techniques.
\end{itemize}

This general technique for the study of the enumerative properties
of pattern avoiding permutations can be properly described in the
framework of the ECO method \cite{BDLPP}, which offers a rigorous
setting for the notions of ``insertion of a node", ``active site"
and ``succession rule" mentioned above.

\bigskip

The above described graphical representation for permutations can
be easily conformed to words. The only difference lies in the fact
that, in the representation of a word, more than one node may
appear on the same line.

\bigskip

In our paper, the general strategy we have briefly sketched above
is suitably extended for our purpose, adapting it to the case of
pattern avoiding words.

\section{Notations and definitions}

Throughout the paper $M$ will denote a finite totally ordered
alphabet, with $|M|=m$. Therefore we will always set $M=\{
1,2,\ldots m\}$.

A \emph{word pattern} (or simply \emph{pattern}) in $M$ is a word
$x$ on an alphabet $\Sigma$, where $\Sigma =\{ 1,\ldots ,k\}$,
with $k\leq m$, where each letter of $\Sigma$ appears in $x$ at
least once. A word pattern in $M$ will also be called a
\emph{reduced word} on $M$. For instance, 32212 is a pattern in
$M=\{ 1,2,\ldots ,10\}$, but 41776 isn't.

Denoting as usual by $M^*$ the set of all the words on the
alphabet $M$, we say that $w\in M^*$ \emph{avoids} a pattern $x$
in $M$ when no subword of $w$ is order-isomorphic to $x$ (see
\cite{B}). For example, $w=311472511\in \{1,\ldots ,7\}^*$ avoids
the patterns 212 and 221, but not the pattern 121.

The notion of generalized pattern for words is analogous to the
definition introduced for permutations; however, the rigorous
definition can be found in \cite{BM}. In particular, a generalized
pattern of length three \emph{of type $(1,2)$} is a reduced word
$x\in \{ 1,2,3\}^*$ of the kind $a-bc$ (for the meaning of dashes
in generalized word patterns, we refer once again to \cite{BM}).

Finally, we will adopt the following notations, which differ from
the corresponding ones used in \cite{B,BM}:
\begin{itemize}
\item $W_n ^{(m)}(T)$: set of words on $M$ of length $n$ avoiding
each pattern in the set $T$; \item $f_T ^{(m)}(x)$: generating
function of $W_n ^{(m)}(T)$ (with respect to the length of the
words).
\end{itemize}

\section{The general method}\label{method}

In this section we give a detailed description of the steps to be
followed in order to apply our method to the enumeration of
classes of (generalized) pattern avoiding words. We remark that
all the notations introduced in the present section will be
extensively used throughout the paper.

Suppose to be interested in counting the words in $W_n ^{(m)}(T)$
according to their length $n$.

\begin{description}
\item \underline{Step 1.} Denote by $\widetilde{W}_n ^{(k)}(T)$
the set of reduced words of $W_n ^{(m)}(T)$ containing each letter
of $K=\{ 1,2,\ldots ,k\}$. Our first goal is to determine
$|\widetilde{W}_n ^{(k)}(T)|$, or, at least, an expression for the
generating function of such quantities. \item \underline{Step 2.}
Now we try to describe an effective construction of the set of the
reduced words of length $n+1$ starting from the set of the reduced
words of length $n$. To do this, we will apply the so-called ECO
method, for which we refer the reader to \cite{BDLPP,DFR,BETAL}.
Taken a word $w\in \widetilde{W}_n ^{(k)}(T)$, we try to insert a
letter to the right of $w$, so obtaining a new word $w'$ of length
$n+1$. This can be done by describing $w$ by means of a path-like
representation completely analogous to the one used for
permutations. Here, the only difference consists of the fact that
we can add the rightmost letter both in a region (if the letter
occurs for the first time in $w$) and on a horizontal line
(otherwise). After the insertion, we possibly need to rename part
of the letters of $w$, according to the added letter. Moreover,
the choice of the letter to be added heavily depends on the
patterns to be excluded. As an example, consider $w=12132\in
\widetilde{W}_5 ^{(3)}(1-22,2-12)$. Using the graphical
representation of words mentioned in the introduction and denoting
with a circle a place where a letter can be added, Figure 1 is
obtained.

\begin{figure}\label{pippone}
\begin{center}
\includegraphics [scale=0.2]{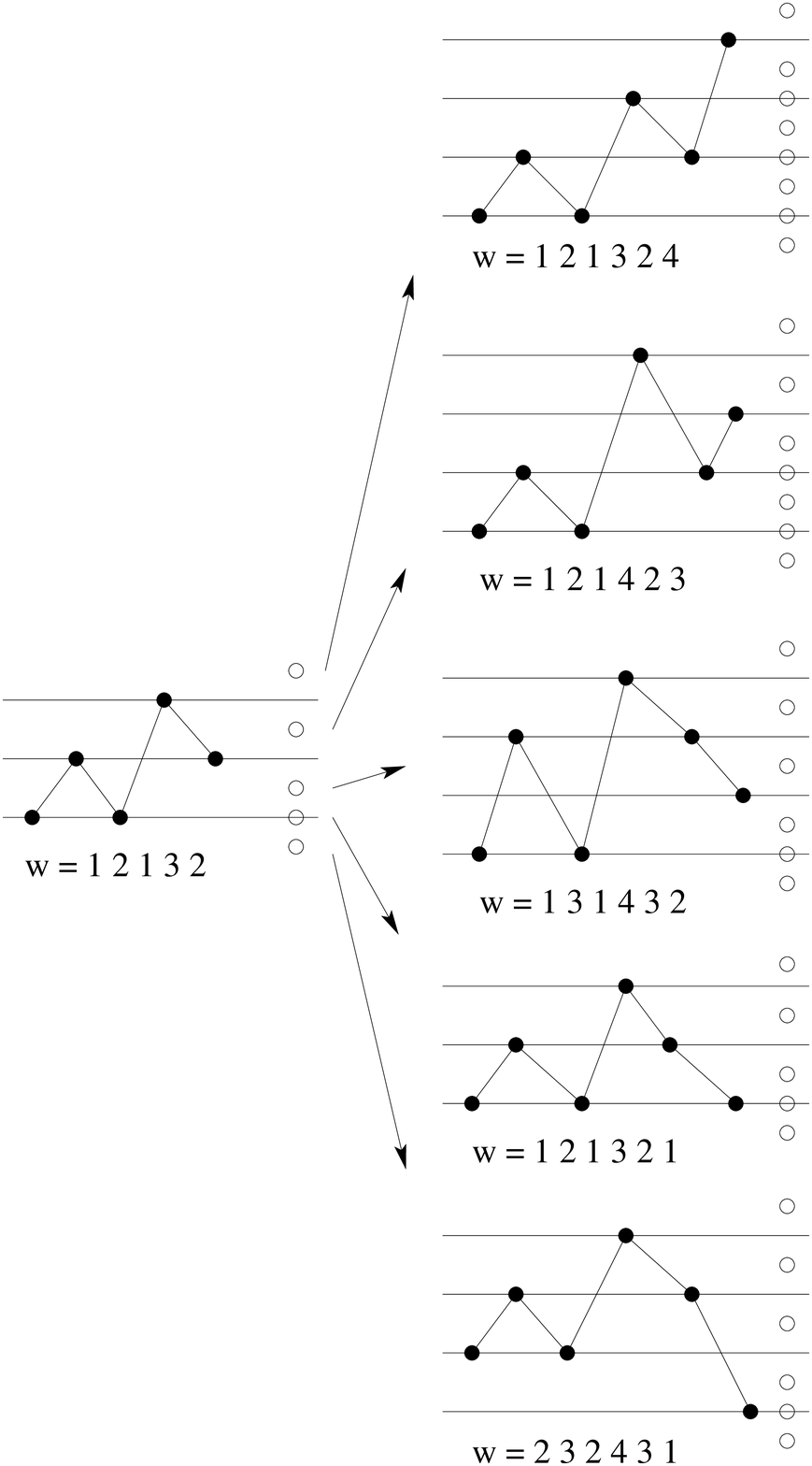}
\end{center}
\caption{the production of $w=12132$}
\end{figure}

Figure 1 means that, starting from 12132, we can construct five
new words (of length 6), which are the ones described in the
picture. As one can immediately notice, if $w\in \widetilde{W}_n
^{(k)}(T)$, then all the words produced by $w$ (which, in the
sequel, will also be called the \emph{sons} of $w$) belong to
$\widetilde{W}_{n+1} ^{(k)}(T)\cup \widetilde{W}_{n+1}
^{(k+1)}(T)$.

\item \underline{Step 3.} If we are lucky, the ECO construction
found in step 2 can be translated into a succession rule or, which
is essentially the same, into a generating tree \cite{W}. In turn,
if the succession rule is regular enough, we can hope to use it to
get a closed formula for the numbers
$\alpha_{n,k}=|\widetilde{W}_n ^{(k)}(T)|$ or maybe, if this is
not possible, to derive a nice expression for the generating
functions $f_k (x)=\sum_n \alpha_{n,k}x^n$, for any fixed $k$.

\item \underline{Step 4.} If we have been able to enumerate
$\widetilde{W}_n ^{(k)}(T)$, for all $k\leq m$, it is now easy to
count the words of $W_n ^{(m)}(T)$. Indeed, it is clear that each
reduced word $w$ in $\widetilde{W}_n ^{(k)}(T)$ is associated with
exactly ${m\choose k}$ distinct words in $W_n ^{(m)}(T)$: just
replace the set of letters of $w$ with any possible subset of $M$
having $k$ elements, taking care of preserving the relative order
of the letters. For instance, the reduced word $1213\in
\widetilde{W}_4 ^{(3)}(1-11,2-11)$ is associated with the
following ${5\choose 3}=10$ words of $W_4 ^{(5)}(1-11,2-11)$:
\begin{eqnarray*}
1213\qquad 1214\qquad 1215\qquad 1314\qquad 1315
\\ 1415\qquad 2324\qquad 2325\qquad 2425\qquad 3435
\end{eqnarray*}

Thus, ${m\choose k}\alpha_{n,k}$ is the number of words of $W_n
^{(m)}(T)$ having $k$ distinct letters. Finally, we have simply to
sum up to get :
\begin{displaymath}
|W_n ^{(m)}(T)|=\sum_{k=0}^{n}{m\choose k}\alpha_{n,k}.
\end{displaymath}
\end{description}

Before applying the above described general method to some
specific classes of pattern avoiding words, an important remark is
to be done. When we construct the generating tree associated with
the ECO construction found in step 2, what we get is not really a
description of the growth of $\widetilde{W}(T)$. Indeed, if
$|M|=m$, the construction depicted in Figure 1 only works when the
starting word belongs to $\widetilde{W}_n ^{(k)}(T)$, {\em for
$k<m$}. If $k=m$, it is clear that, in our graphical
representation, no new horizontal line can be added in order to
generate new words. For instance, for $m=3$ and $T=\{
1-22,2-12\}$, the only son of the word 12132 is the word 121321
(see figure below).

\begin{center}
\includegraphics [scale=0.3]{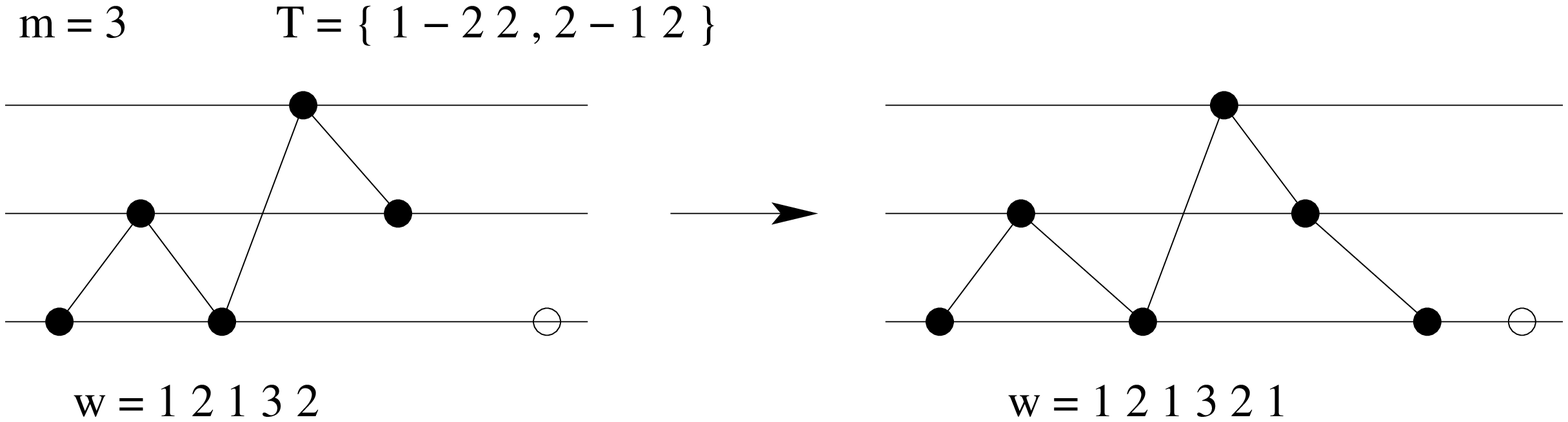}
\end{center}

Therefore, every time we will apply our construction, we have to
keep in mind that, for $k<m$, it is possible to add new horizontal
lines, whereas for $k\geq m$ no further horizontal line can be
added: hence our generation procedure, and the strictly related
enumeration technique, must be suitably adapted.

\section{The class $W(1-12,2-21)$}

In this section we will apply our method to enumerate the pattern
avoiding class $W_n ^{(m)} (1-12,2-21)$. As we will see, this will
lead us to consider a succession rule having only odd labels and,
in particular, an infinite number of labels producing only one
son.

Following our program, we start by considering the set
$\widetilde{W}_n ^{(k)} (1-12,2-21)$ of the reduced words of
$\widetilde{W}_n ^{(m)} (1-12,2-21)$ on the $k$-letter alphabet
$\{ 1,2,\ldots ,k\}$. Our first goal is to determine
$|\widetilde{W}_n ^{(k)} (1-12,2-21)|$.

Take $w\in \widetilde{W}_n ^{(k)} (1-12,2-21)$ and suppose that
the last letter of $w$ is $h$ ($\leq k$). We can distinguish two
cases:
\begin{itemize}
\item[i)] The letter $h$ occurs in $w$ for the first time in the
last position. We represent this situation using the following
diagram:

\begin{center}
\includegraphics [scale=0.3]{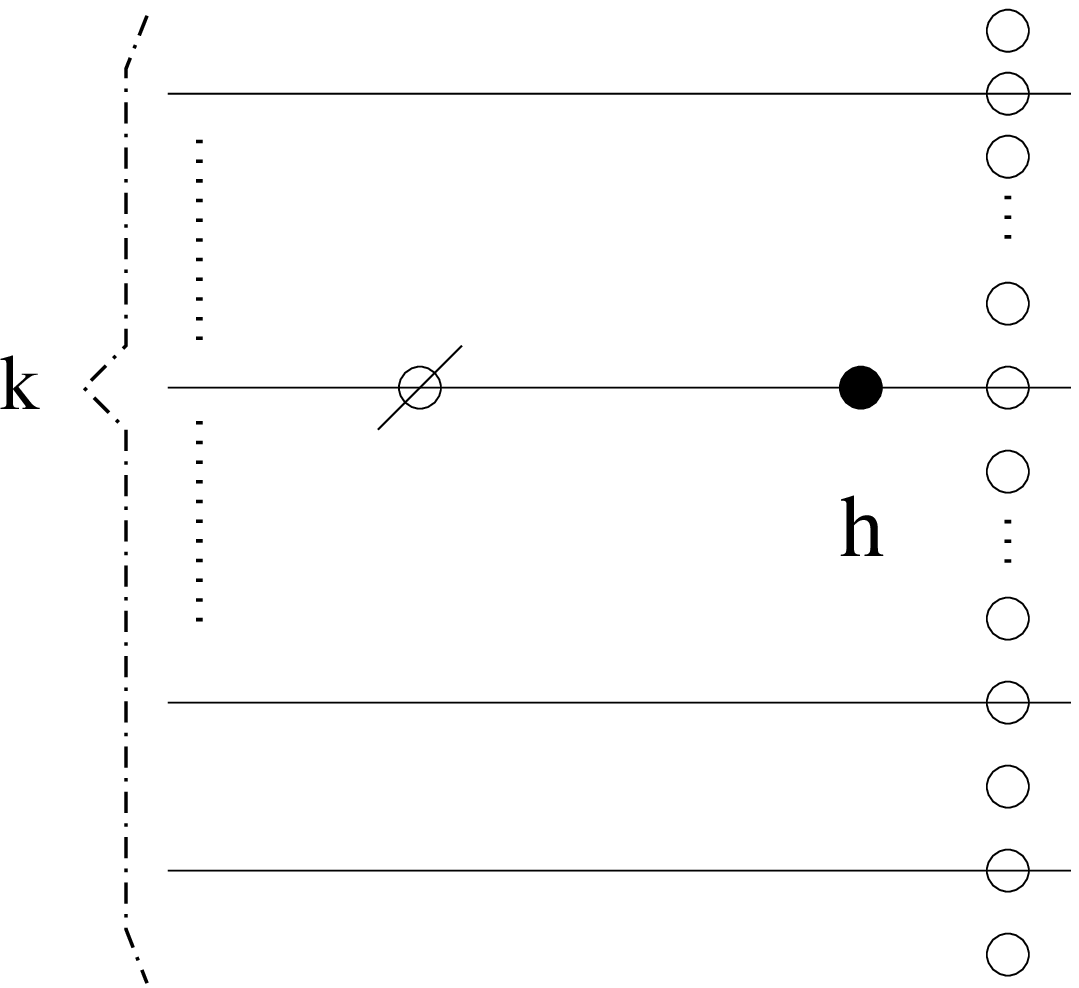}
\end{center}

The above diagram should be interpreted as follows: $w$ is a
reduced word on a $k$-letter alphabet; the empty set symbol in the
same line of the last letter denote that the $h$ in the last
position is the first occurrence of $h$ in $w$; finally, the
circles on the right indicate the active sites where we can add a
new letter in order to obtain a word $w'\in \widetilde{W}_{n+1}
^{(k)} (1-12,2-21)\cup \widetilde{W}_{n+1} ^{(k+1)} (1-12,2-21)$.
In this particular case, we are allowed to add any letter of $\{
1,2,\ldots ,k+1\}$ on the right of $w$. To translate this fact
into a succession rule, we use the following notations: the label
$(k_h)$ denotes a word of $\widetilde{W}_n ^{(k)} (1-12,2-21)$
ending with the letter $h$ and having no further occurrences of
$h$ before; the label $(\overline{k}_h)$ denotes a word of
$\widetilde{W}_n ^{(k)} (1-12,2-21)$ ending with the letter $h$
and such that $h$ also appears in some previous position. With
these notations, the production rule encoding the above
construction can be written as follows:
\begin{displaymath}(k_h)\rightsquigarrow (\overline{k}_1)\cdots (\overline{k}_k)((k+1)_1)\cdots
((k+1)_{k+1}).
\end{displaymath}

\item[ii)] The letter $h$ also occurs in $w$ in some position
other than the last one. In this case, we can use the following
diagram:

\begin{center}
\includegraphics [scale=0.3]{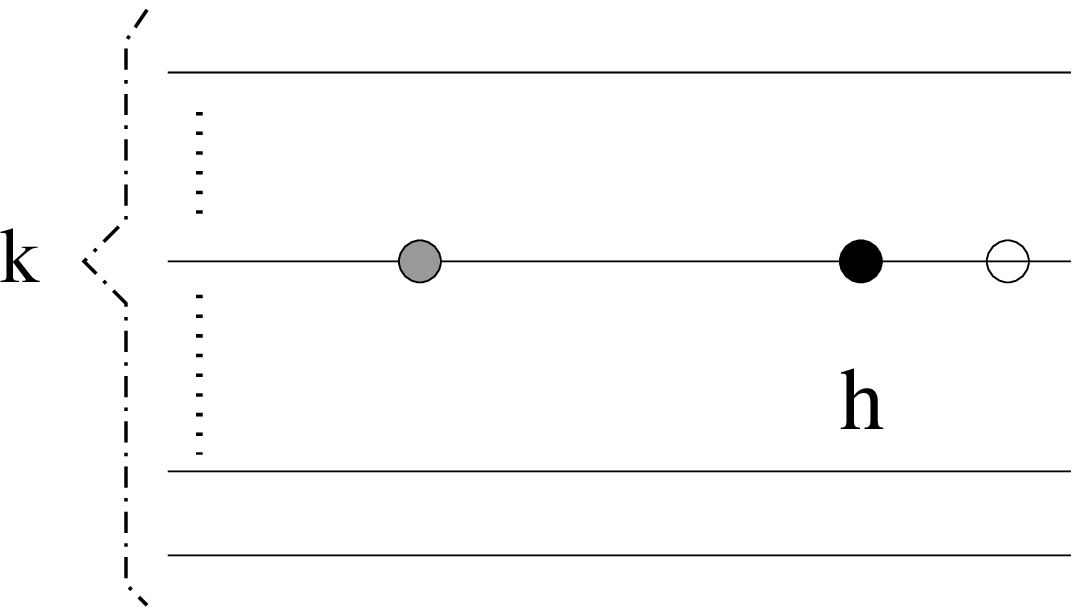}
\end{center}

The meaning of the symbols are obviously the same as for the
previous diagram. Here the only difference consists of the fact
that $h$ appears not only in the last position of $w$, but also
somewhere before, and this has been represented by a gray circle
placed in the same line of the last occurrence of $h$. In this
second case, every insertion in a place (region or line) other
than such a line would produce an occurrence of a forbidden
pattern. Therefore, with the same notations as above, the
associated production rule is the following:
\begin{displaymath}
(\overline{k}_h)\rightsquigarrow (\overline{k}_h).
\end{displaymath}
\end{itemize}

Now, putting things together, and recalling that the word of
minimum length avoiding simultaneously $1-12$ and $2-21$ is 1
(which is represented by the label $(1_1)$), the succession rule
encoding the whole construction is the following:
\begin{eqnarray}\label{complic}
\left\{ \begin{array}{lll} (1_1)
\\ (k_h)\rightsquigarrow (\overline{k}_1)\cdots (\overline{k}_k)((k+1)_1)\cdots ((k+1)_{k+1})
\\ (\overline{k}_h)\rightsquigarrow (\overline{k}_h)
\end{array}\right. .
\end{eqnarray}

The above succession rule is indeed too complicated to allow
enumeration. However, it is possible to rewrite it in an
equivalent form which is surely more suitable for our purposes. To
do this, the following two observations are essential.
\begin{enumerate}
\item In the production $(\overline{k}_h)\rightsquigarrow
(\overline{k}_h)$, the subscript $h$ is unimportant; instead, we
should keep track of the parameter $k$, since we will use it in
the enumeration process. So, since $(\overline{k}_h)$ has only one
production, we can replace this production rule with the rule
$(1_k)\rightsquigarrow (1_k)$ (so that $(1_k)$ stands for a
reduced word on a $k$-letter alphabet whose last letter appears
somewhere else in the word itself). \item The production of
$(k_h)$ is independent of $h$; therefore, since $(k_h)$ produces
$2k+1$ labels, taking care of the previous remark, we can replace
its production rule with the following:
\begin{displaymath}
(2k+1)\rightsquigarrow (1_k)^k (2k+3)^{k+1}.
\end{displaymath}

Thus, $(2k+1)$ represents a $k$-letter word on a reduced alphabet
whose last letter does not appear in any other position.
\end{enumerate}

Thanks to the above considerations, we have that the rule in
(\ref{complic}) is equivalent to:
\begin{equation}\label{sempl}
\left\{ \begin{array}{lll} (3)
\\ (2k+1)\rightsquigarrow (1_k)^k (2k+3)^{k+1}
\\ (1_k)\rightsquigarrow (1_k)
\end{array}\right. .
\end{equation}

The first levels of the associated generating tree are depicted in
the figure below:

\begin{center}
\includegraphics [scale=0.4]{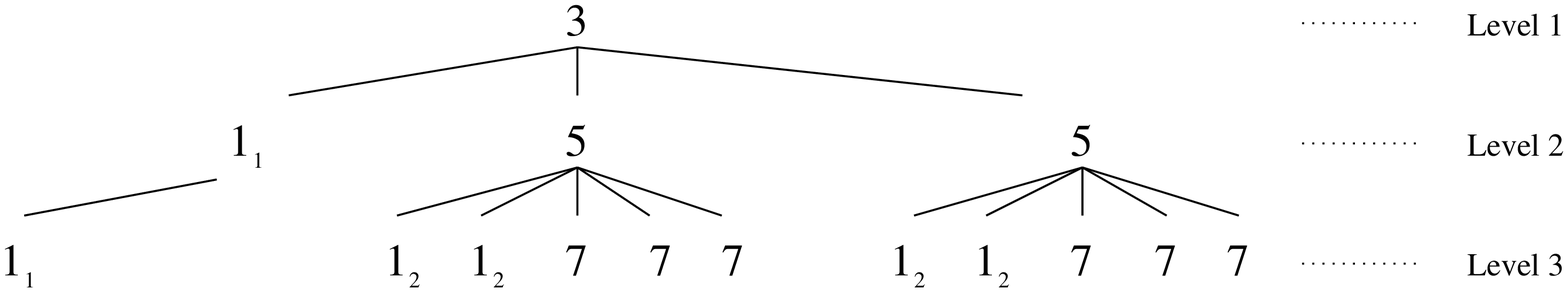}
\end{center}

Looking at the labels of the succession rule (\ref{sempl}) and
keeping in mind how such a rule has been obtained, we immediately
observe that:
\begin{itemize}
\item[i)] each label $(2k+1)$ represents a word having $k$
letters; \item[ii)] each label $(1_k)$ represents a word having
$k$ letters; \item[iii)] the unique labels appearing at level $n$
of the associated generating tree are $(1_1),(1_2),\ldots
,(1_{n-1}),(2n+1)$.
\end{itemize}

Therefore, $\alpha_{n,k}$ (see Section 3) is either the number of
labels $(1_k)$ at level $n$, if $k<n$, or the number of labels
$(2k+1)$ at level $n$, if $k=n$.

From the succession rule (\ref{sempl}), we can obtain an explicit
formula for the number $\alpha_{n,k}$.

\begin{teor} When $k\leq m$, we have
\begin{displaymath}
\alpha_{n,k}=\left\{ \begin{array}{ll} k\cdot k!\; ,&k<n
\\ n!\; ,&k=n
\end{array}\right. ,
\end{displaymath}
whereas, for $k>m$, it is $\alpha_{n,k}=0$.
\end{teor}

\emph{Proof.} We can use the rule (\ref{sempl}) to get the
following recursions for the $\alpha_{n,k}$'s:
\begin{eqnarray}
&&\alpha_{n,n}=n\cdot \alpha_{n-1,n-1}\label{uno}
\\ &&\alpha_{n,n-1}=(n-1)\cdot \alpha_{n-1,n-1}\label{due}
\\ &&\alpha_{n,k}=\alpha_{k+1,k}\qquad (k<n-1) .\label{tre}
\end{eqnarray}

Since $\alpha_{1,1}=1$, from (\ref{uno}) we immediately have
$\alpha_{n,n}=n!$, whence, from (\ref{due}), we get
$\alpha_{n,n-1}=(n-1)\cdot (n-1)!$. Finally, from (\ref{tre}), it
is $\alpha_{n,k}=k\cdot k!$, for $k<n-1$. These formulas hold only
for $k\leq m$, since our alphabet has $m$ letters, and so the
generating tree of (\ref{sempl}) can be used only for $k<m$.
Obviously, when $k>m$ we have $\alpha_{n,k}=0$.\cvd

The first lines of the infinite matrix $A=(\alpha_{n,k})_{n,k\geq
0}$ are the following (where we have added one more row and one
more column to represent the empty word):
\begin{displaymath}
A=\begin{pmatrix} 1&0&0&0&0&0&\cdots
\\ 0&1&0&0&0&0&\cdots
\\ 0&1&2&0&0&0&\cdots
\\ 0&1&4&6&0&0&\cdots
\\ 0&1&4&18&24&0&\cdots
\\ 0&1&4&18&96&120&\cdots
\\ \vdots&\vdots&\vdots&\vdots&\vdots&\vdots&\ddots
\end{pmatrix}.
\end{displaymath}

Of course, in the matrix $A$ we have to turn into $0$ the entries
of all the columns from the $(m+1)$-th onwards.

\bigskip

Now we are ready to get the exact enumeration of the class
$W(1-12,2-21)$. Following the method described in section
\ref{method}, we have the following results.

\begin{prop} The total number of words of $W_n ^{(m)}(1-12,2-21)$ having $k\leq m$ distinct letters
is:
\begin{displaymath}
{m\choose k}\alpha_{n,k}=\left\{ \begin{array}{ll} k\cdot (m)_k \;
,&k<n
\\ (m)_n \; ,&k=n
\end{array}\right. ,
\end{displaymath}
where $(a)_b$ denotes the usual falling factorial, $(a)_b =a\cdot
(a-1)\cdot \ldots \cdot (a-b+1)$.
\end{prop}

\begin{teor} It is
\begin{displaymath}
|W_n ^{(m)} (1-12,2-21)|=\left\{ \begin{array}{ll}
\sum_{k=0}^{n-1}k\cdot (m)_k +(m)_n \; ,& n\leq m
\\ \sum_{k=0}^{m-1}k\cdot (m)_k \; ,& n>m
\end{array}\right. .
\end{displaymath}
\end{teor}

\section{The class $W(1-21,2-12)$}

The next case we take into consideration is that of words
simultaneously avoiding the two generalized patterns $1-21$ and
$2-12$. Also in this case, our technique allows us to get to an
explicit counting result quite quickly and easily.

As usual, we start by considering $\widetilde{W}_n
^{(k)}(1-21,2-12)$, i.e. reduced words. Suppose that $w\in
\widetilde{W}_n ^{(k)}(1-21,2-12)$ and that $h$ ($\leq k$) is the
last letter of $w$. Using a diagram similar to the one of the
previous section, we can represent the set of productions of $w$
as follows:

\begin{center}
\includegraphics [scale=0.3]{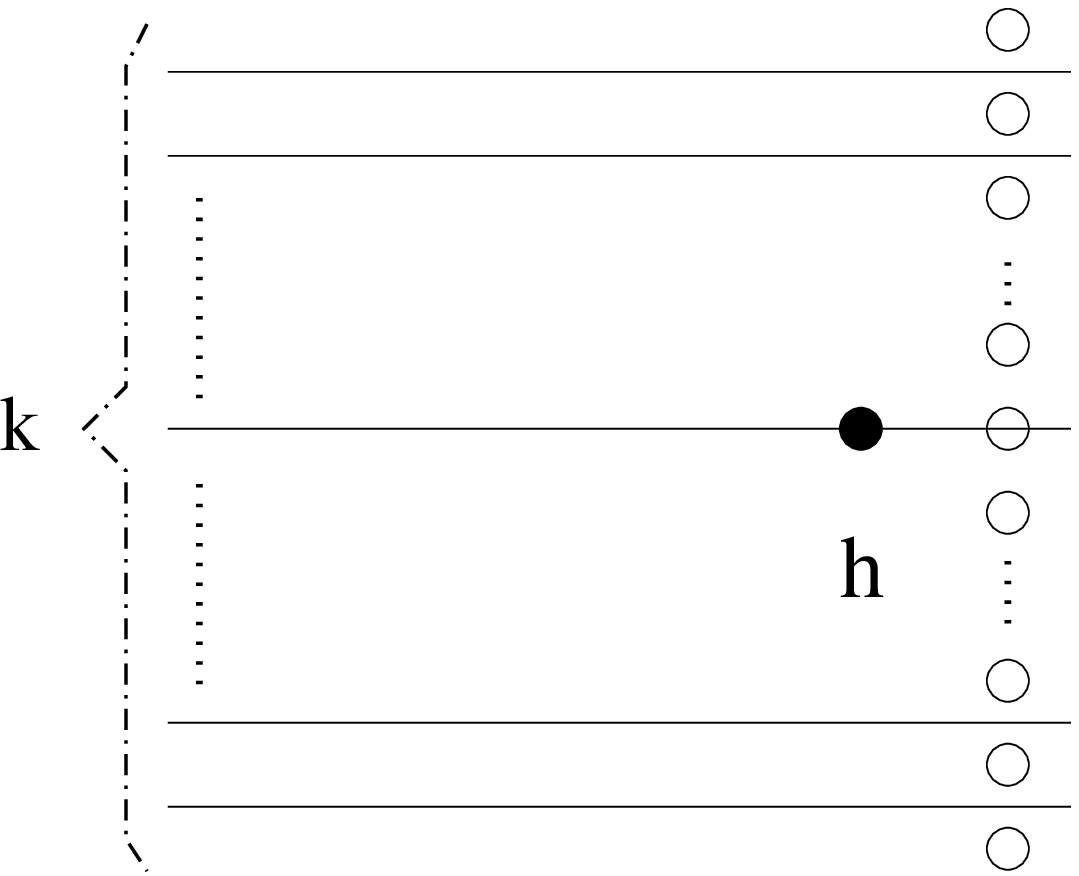}
\end{center}

Notice that the insertion is not possible on a line other than the
one containing the last letter of $w$.

Unlike the previous class of words, here it is not necessary to
distinguish two cases, since the production of $w$ is independent
of the fact that the last letter of $w$ is the first occurrence of
$h$ or not. Since $1\in \widetilde{W}_n ^{(k)}(1-21,2-12)$ is the
only word of length 1 of our class, we have the following
succession rule associated with the above diagram:
\begin{displaymath}
\left\{ \begin{array}{ll} (1_1 )
\\ (k_h )\rightsquigarrow (k_h )((k+1)_1 )\cdots ((k+1)_{k+1})
\end{array}\right. .
\end{displaymath}

Also in this case, we can observe that the number of labels
produced by $(k_h )$ does not depend on $h$. Since we are only
interested in the parameter $k$, we can easily find that the above
rule is equivalent to
\begin{equation}\label{semplice}
\left\{ \begin{array}{ll} (3)
\\ (k)\rightsquigarrow (k)(k+1)^{k-1}
\end{array}\right. .
\end{equation}

In the rule (\ref{semplice}), a label $(k)$ represents a word
having $k-2$ distinct letters. The first few levels of the
associated generating tree look as follows:

\begin{center}
\includegraphics [scale=0.4]{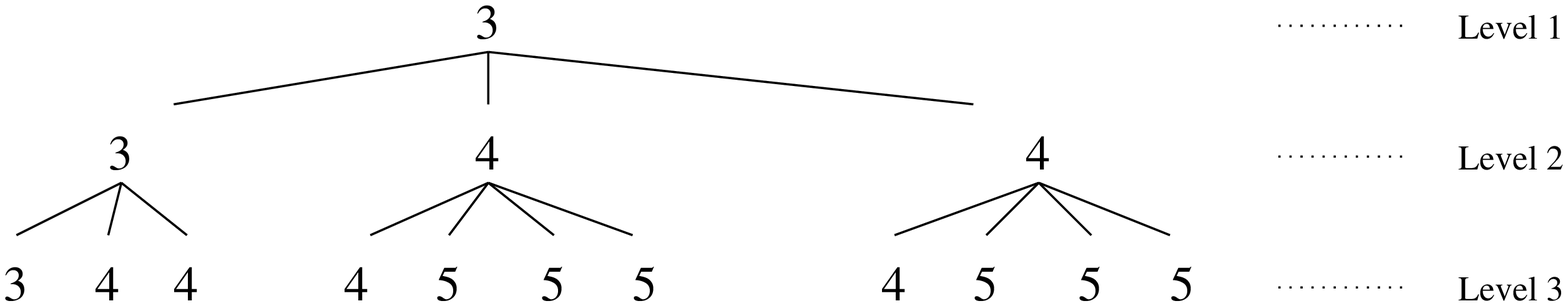}
\end{center}

The above succession rule allows us to find the exact value of the numbers
$\alpha_{n,k}=|\widetilde{W}_n ^{(k)}(1-21,2-12)|$.

\begin{prop} For $k\leq m$, we have
\begin{equation}
\alpha_{n,k}=k\cdot (n-1)_{k-1}.
\end{equation}
\end{prop}

\emph{Proof.} The same production rule as (\ref{semplice}), but
with axiom 2, constitutes a well known succession rule, where the
associated numerical sequence is given by the number of
arrangements \cite{BETAL}. It is immediate to realize that our
generating tree can be seen as a part of the generating tree for
arrangements, as the following figure clarifies:

\begin{center}
\includegraphics [scale=0.4]{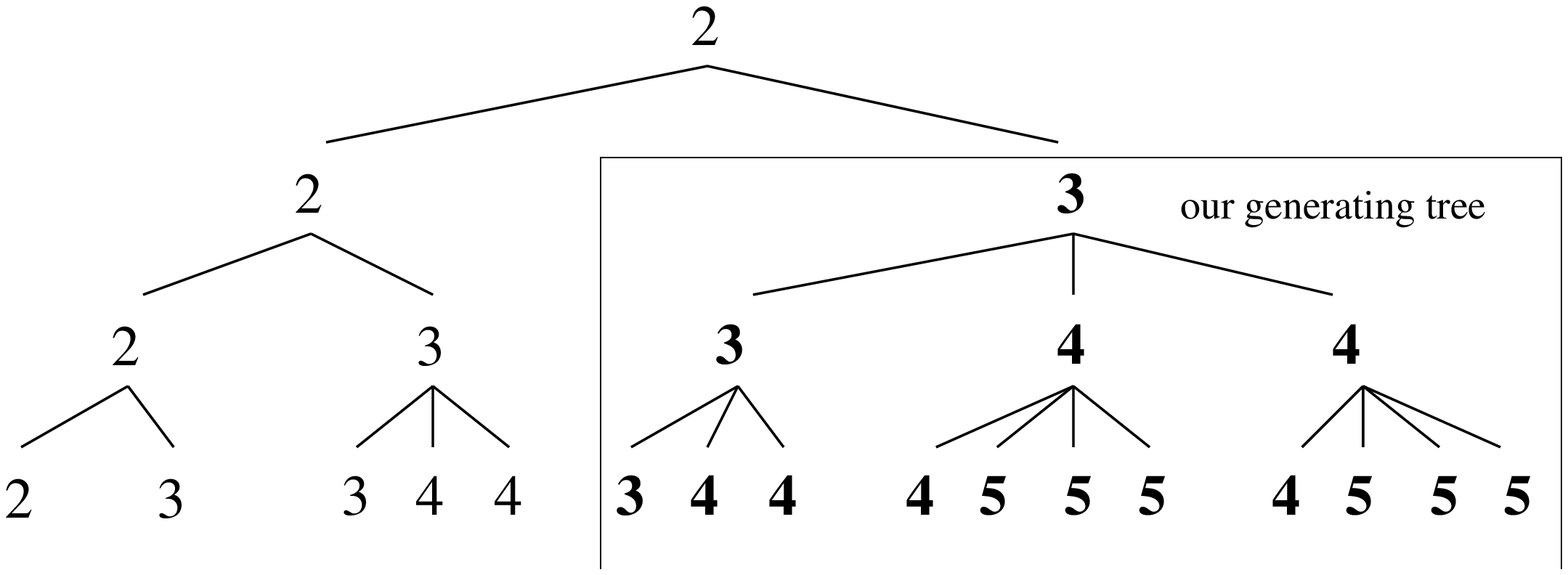}
\end{center}

The level polynomials of the generating tree for arrangements are
\begin{displaymath}
p_n (x)=\sum_{k\geq 0}(n-1)_k x^k
\end{displaymath}
(here the coefficient of $x^k$, depending on $n$, gives the number
of labels $k+2$ at level $n$). Therefore, the relation between
$p_n (x)$ and the level polynomials $\alpha_n (x)$ (where the
coefficient of $x^k$ is the number of labels $k+3$ at level $n$)
of our generating tree is expressed by the equality
\begin{displaymath}
p_n (x)=p_{n-1}(x)+x\alpha_{n-1}(x),
\end{displaymath}
which gives
\begin{displaymath}
\alpha_n (x)=\frac{p_{n+1}(x)-p_n (x)}{x},
\end{displaymath}
whence, recalling that the coefficients of the polynomials
$\alpha_n (x)$ coincide with the numbers $\alpha_{n,k}$ only when
$k\leq m$, for such values of $k$ we have
\begin{displaymath}
\alpha_{n,k}=(n)_k -(n-1)_k =k\cdot (n-1)_{k-1}.\quad \blacksquare
\end{displaymath}

\bigskip

The first few lines of the matrix $A=(\alpha_{n,k})_{n,k\geq 0}$
(where the first column, representing the set of words having $0$
letters, and the first row have been added) are:
\begin{displaymath}
A=\begin{pmatrix} 1&0&0&0&0&0&\cdots
\\ 0&1&0&0&0&0&\cdots
\\ 0&1&2&0&0&0&\cdots
\\ 0&1&4&6&0&0&\cdots
\\ 0&1&6&18&24&0&\cdots
\\ 0&1&8&36&96&120&\cdots
\\ \vdots&\vdots&\vdots&\vdots&\vdots&\vdots&\ddots
\end{pmatrix}.
\end{displaymath}

As usual, we must remember that the above formula is only valid
when $k\leq m$, so that in the matrix $A$ we have to consider
columns exclusively up to $k=m$, setting all the remaining entries
equal to $0$.

Finally we can give complete enumeration results.

\begin{prop} The number of words of $W_n ^{(m)}(1-21,2-12)$ having $k$ distinct
letters is
\begin{displaymath}
{m\choose k}\alpha_{n,k}={m\choose k}k\cdot (n-1)_{k-1}.
\end{displaymath}
\end{prop}

\begin{teor} The total number of words of $W_n ^{(m)}(1-21,2-12)$ is
\begin{displaymath}
|\widetilde{W}_n ^{(k)}(1-21,2-12)|=\sum_{k=0}^{n}{m\choose
k}k\cdot(n-1)_{k-1}.
\end{displaymath}
\end{teor}

\section{The class $W(1-11,1-12)$}

This case is more difficult to deal with than the previous ones.
Of course, we will follow the same general method; at the end, we
will get an expression for the generating function of the class,
rather than a closed form for the associated numerical sequence.

As usual, we start by considering the set $\widetilde{W}_n
^{(k)}(1-11,1-12)$ of $k$-reduced words of length $n$
simultaneously avoiding $1-11$ and $1-12$. Let $w\in
\widetilde{W}_n ^{(k)}(1-11,1-12)$ with last letter $h$. We can
distinguish two cases.
\begin{itemize}
\item[i)] The occurrence of $h$ in the last position of $w$ is the
first occurrence of $h$ in $w$. Drawing the usual diagram, we are
in the following situation:

\begin{center}
\includegraphics [scale=0.3]{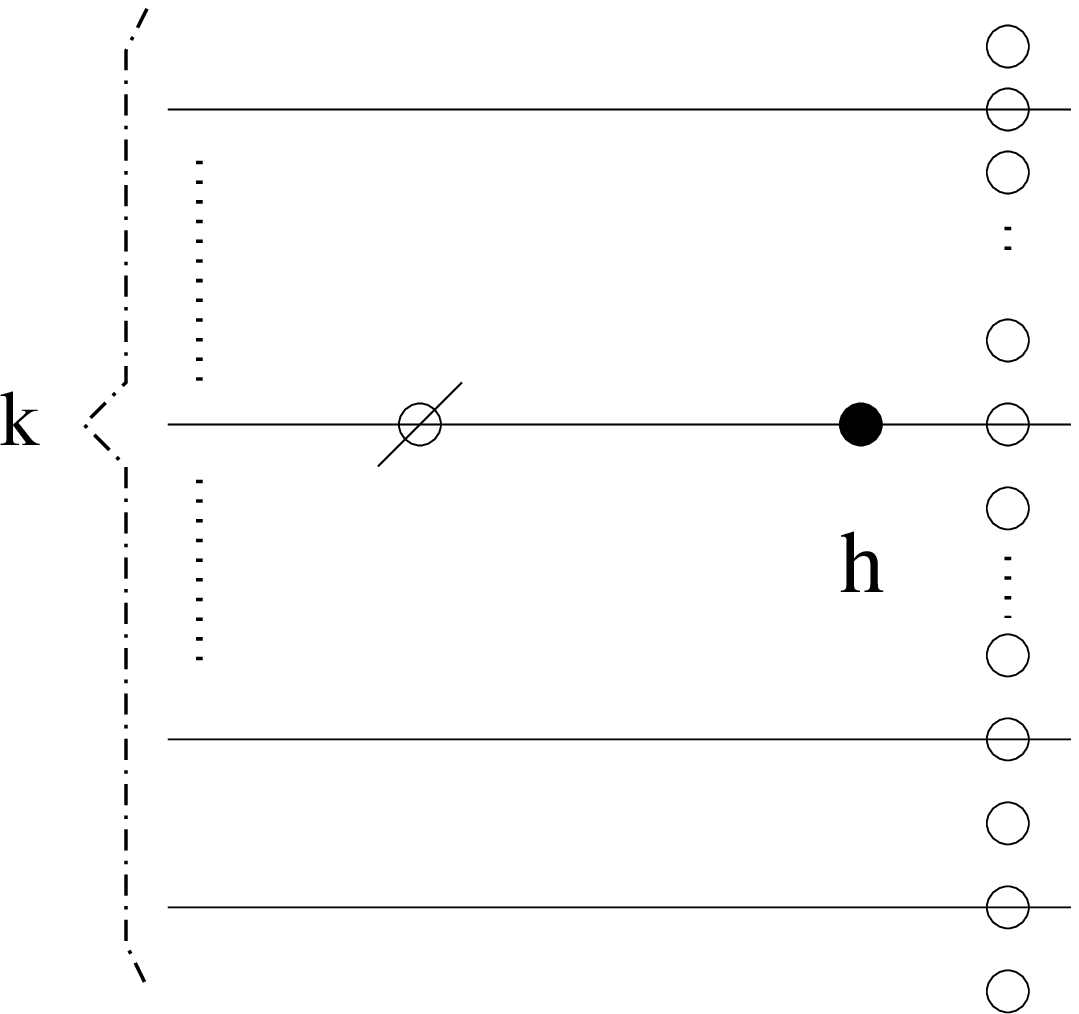}
\end{center}

Remember that the circles on the right represent those active
sites in which we are allowed to add a letter at the end of $w$.
In this case, the above diagram tells us that we can add a letter
in any position. Thus we can encode this fact using the following
production rule:
\begin{equation}\label{prod1}
(k_h )\rightsquigarrow (\overline{k}_1 )\cdots (\overline{k}_k
)((k+1)_1 )\cdots ((k+1)_{k+1})
\end{equation}
(the meaning of the symbols has been extensively explained in the
previous two cases).

\item[ii)] The letter $h$ occurs in $w$ for the first time
somewhere before its occurrence in the last position. In this case
the situation is the following:

\begin{center}
\includegraphics [scale=0.3]{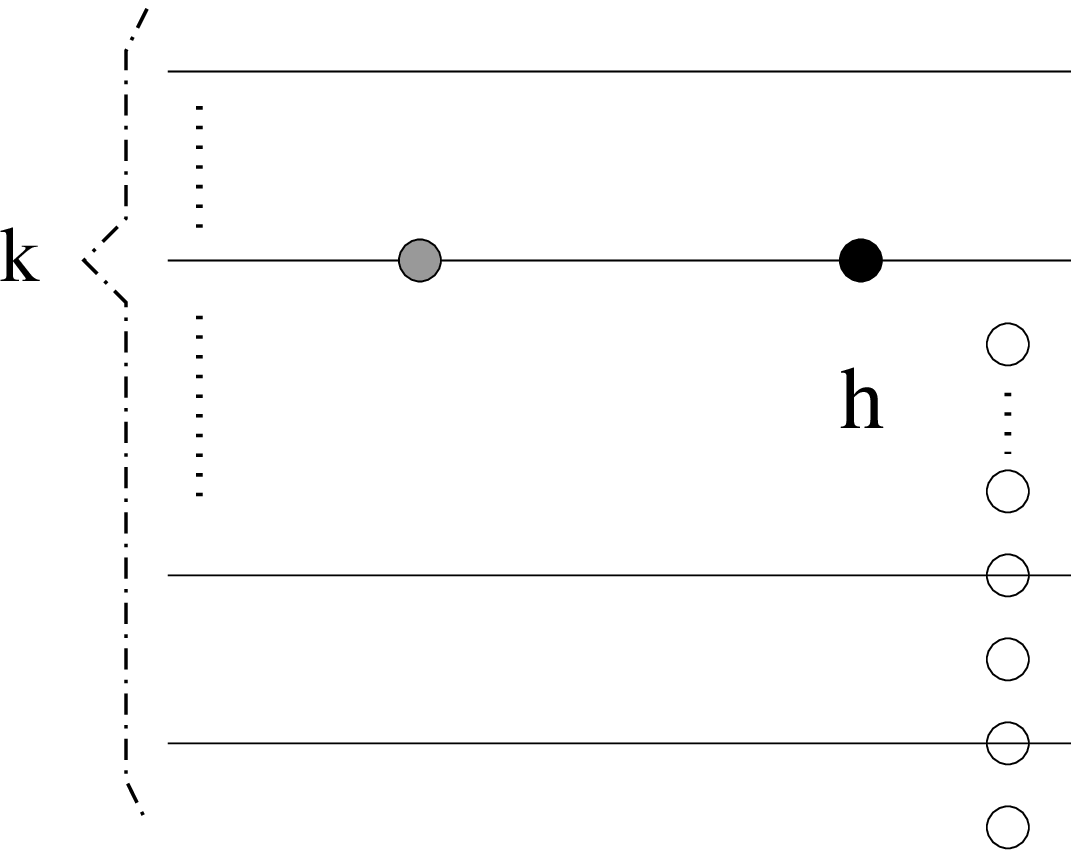}
\end{center}

Here the production is not as trivial as in case i), and can be
encoded by the rule
\begin{equation}\label{prod2}
(\overline{k}_h)\rightsquigarrow (\overline{k}_1)\cdots
(\overline{k}_{h-1})((k+1)_1 )\cdots ((k+1)_h ).
\end{equation}
\end{itemize}

What we can immediately observe is that, in case i), the
production of $(k_h )$ does not depend on $h$, whereas, in case
ii), the production of $(\overline{k}_h )$ does. Moreover, setting
$h=k+1$ both in (\ref{prod1}) and in (\ref{prod2}), we obtain
formally the same production rule. This implies that our ECO
construction for the class $\widetilde{W}_n ^{(k)}(1-11,1-12)$ is
equivalent to the following:
\begin{equation}\label{regola}
\left\{ \begin{array}{ll} (1_2 )
\\ (k_h )\rightsquigarrow (k_1 )\cdots (k_{h-1})((k+1)_{k+2})^{h}\qquad ,\qquad h\leq k+1
\end{array}\right. ,
\end{equation}
where the label $(k_h )$, for $h<k+1$, corresponds to the labels
$(\overline{k}_h )$ in (\ref{prod1}) and (\ref{prod2}) and the
label $(k_{k+1})$ corresponds to the label $(k_h) $ in
(\ref{prod1}) and (\ref{prod2}), {\em for any $h\leq k+1$}. We can
also have a look at the first levels of the associated generating
tree:

\begin{center}
\includegraphics [scale=0.4]{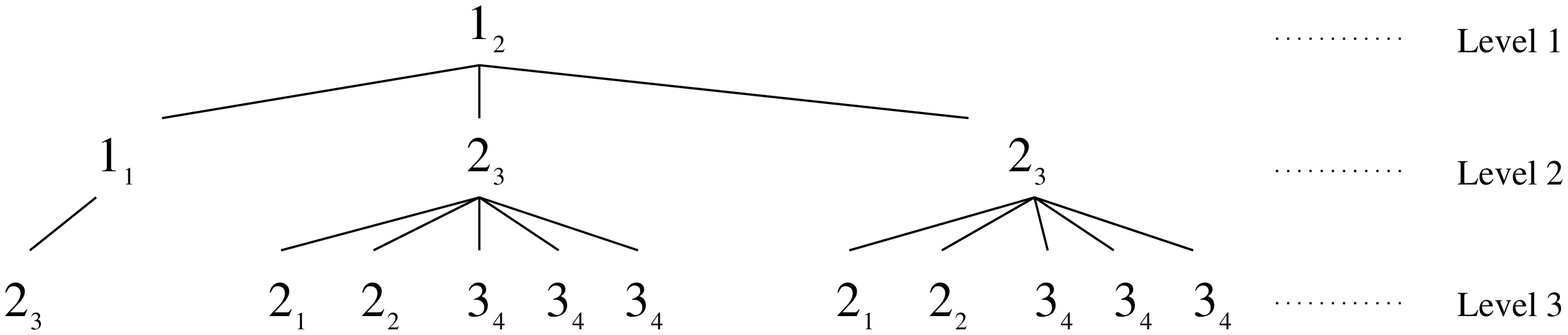}
\end{center}

Using the succession rule (\ref{regola}), we are able to enumerate
the class $\widetilde{W}_n ^{(k)}(1-11,1-12)$. To do this, we will
make use of a recursion on the entries of the associated ECO
matrix we are going to show below.

The first lines of the ECO matrix associated with the above
depicted generating tree (i.e., the infinite matrix describing the
distribution of the labels at the various levels of the generating
tree, see \cite{DFR}) can be organized as follows:

\bigskip

\scriptsize
\begin{center}
\begin{tabular}{|c|cc|ccc|cccc|ccccc|cc|}
\hline

\tiny{labels}$\rightarrow$& $1_2$ & $1_1$ & $2_3$ & $2_2$ & $2_1$
& $3_4$ & $3_3$ & $3_2$ & $3_1$ &
 $4_5$ & $4_4$ & $4_3$ & $4_2$ & $4_1$ & $5_6$ & $\ldots$\\
\tiny{at level} $\downarrow$ &&&&&&&&&&&&&&&&\\

\hline

&&&&&&&&&&&&&&&&\\
1&1&0&&&&&&&&&&&&&&\\
2&0&1&2&0&0&&&&&&&&&&&\\
3&0&0&1&2&2&6&0&0&0&&&&&&&\\
4&0&0&0&1&3&9&6&6&6&24&0&0&0&0&&\\
5&0&0&0&0&1&5&9&15&21&72&24&24&24&24&120&$\cdots$\\
\vdots&$\vdots$&$\vdots$&$\vdots$&$\vdots$&$\vdots$&$\vdots$&$\vdots$&$\vdots$&$\vdots$&
$\vdots$&$\vdots$&$\vdots$&$\vdots$&$\vdots$&$\vdots$&$\vdots$\\
\hline
\end{tabular}
\end{center}
\normalsize

\bigskip

As it should be clear, since $(k_h)$ represents a word having $k$
distinct letters, the cardinality of the set of words having $k$
distinct letters is given by the sum of the entries of the columns
labelled $k_h$, for $h$ running from 1 to $k+1$. To get the
desired recurrence relation among the entries of the matrix, we
recall that the columns labelled $k_h$ ($1\leq h\leq k$) derive
from the production rule (\ref{prod2}), whereas those labelled
$k_{k+1}$ derive from the production rule (\ref{prod1}). In the
sequel  we will denote by $C_k (x)$ the generating function of
column $k_{k+1}$ and by $C_k ^{(h)}(x)$ the generating function of
column $k_h$, with $h\neq k+1$. Columns labelled $k_{k+1}$ will
also be called \emph{special columns}. To keep track of the empty
word we will also consider the generating function $C_0(x)=1$.

\begin{prop} The following recurrence relation holds:
\begin{eqnarray*}
  C_0(x) &=& 1 \\
  C_k(x) &=& \sum_{j=1}^{k}{k\choose j}x^j C_{k-1}(x)\ , \quad k\geq1 \ .
\end{eqnarray*}
\end{prop}

\emph{Proof.} Looking at the succession rule (\ref{regola}), we
observe that a label $(k_h)$, for $1\leq h\leq k$, can only be
generated by a label $(k_t)$, with $t>h$, and each of these labels
produces $(k_h)$ precisely once. Moreover, a label $(k_{k+1})$ can
be generated by any label $((k-1)_h)$, and $((k-1)_h)$ produces
precisely $h$ copies of $(k_{k+1})$. Starting from these
considerations, we are led to the following set of recursions for
the entries of our ECO matrix:
\begin{eqnarray*}
a(1,1_2 )&=&1,\qquad a(n,1_2 )=0\quad (n>1),
\\ a(n,k_h )&=&\sum_{i=h+1}^{k+1}a(n-1,k_i )\quad (h\leq k),
\\ a(n,k_{k+1})&=&\sum_{i=1}^{k}a(n-1,(k-1)_i
)+\sum_{i=1}^{k-1}a(n,(k-1)_i ),
\end{eqnarray*}
where, of course, $a(n,k_h )$ denotes the entry corresponding to
row $n$ and column labelled $k_h$.

The recursion for $a(n,k_h )$, with $h\leq k$, can be iterated so
to obtain a formula expressing $a(n,k_h )$ only in terms of the
entries of special columns, which is:
\begin{eqnarray}\label{ricorsione}
a(n,k_h )=\sum_{j=0}^{k-h}{k-h\choose j}a(n-1-j,k_{k+1}).
\end{eqnarray}
The above formula can be proved using an easy induction argument
on $h$.

\medskip

Now we can plug the above expression for $a(n,k_h )$ into the
recursion formula for $a(n,k_{k+1})$, thus obtaining:
\begin{eqnarray*}
a(n,k_{k+1})&=&a(n-1,(k-1)_k )
\\ &+&\sum_{i=1}^{k-1}\sum_{j=0}^{k-1-i}{k-1-i\choose
j}\left( a(n-2-j,(k-1)_k )+a(n-1-j,(k-1)_k \right).
\end{eqnarray*}

The above recursion for the entries of special columns can be
immediately translated into the desired recursion for the
generating function $C_k (x)$.\cvd

\begin{cor}
\begin{equation}\label{colspecespl}
C_k (x)=\sum_{i=k}^{\frac{k^2 +k}{2}}\left(\sum_{j_1,\cdots,j_k>0
\atop j_1+\cdots +j_k=i}{k\choose j_k }{k-1\choose j_{k-1}}\cdot
\ldots \cdot {1\choose j_1 }\right) x^i .
\end{equation}
\end{cor}

Formula (\ref{colspecespl}) is a completely explicit expression
for the generating function $C_k (x)$, but it is of course
impossible to use it in any expression in which $C_k (x)$ is
required. For this reason, in the sequel we will simply write $C_k
(x)$ in our formulas, but the reader should remember that $C_k
(x)$ can be replaced by its expression in formula
(\ref{colspecespl}).

In order to completely know the entries of our ECO matrix, we now
need to find an expression for the generating functions $C_k
^{(h)}(x)$.

\begin{prop} The following formula (expressing $C_k ^{(h)}(x)$ in
terms of $C_k (x)$) holds:
\begin{displaymath}
C_k ^{(h)}(x)=x(1+x)^{k-h}C_k (x).
\end{displaymath}
\end{prop}

\emph{Proof.} Converting into generating functions the recurrence
relation (\ref{ricorsione}) immediately gives the required
formula.\cvd

Now we are ready to reach our first goal, that is the enumeration
of the class of words $\widetilde{W}_n ^{(k)}(1-11,1-12)$.

\begin{teor} The generating function $f_k (x)$ of
$k$-reduced words of $W_n ^{(m)}(1-11,1-12)$ according to the
length is
\begin{displaymath}
f_k (x)=(1+x)^k C_k (x).
\end{displaymath}
\end{teor}

\emph{Proof.} Since $k$-reduced words of $W_n ^{(m)}(1-11,1-12)$
are represented by the labels of the type $(k_h )$, with $1\leq
h\leq k+1$, the generating function $f_k (x)$ can be obtained by
simply summing up the generating functions corresponding to such
labels, whence
\begin{eqnarray*}
f_k (x)&=&C_k (x)+\sum_{h=1}^k C_k ^{(h)}(x)
\\ &=&C_k (x)+\sum_{h=1}^k x(1+x)^{k-h}C_k (x)
\\ &=&C_k (x)\cdot \left( 1+x\cdot \sum_{h=0}^{k-1}(1+x)^h \right)
\\ &=&(1+x)^k C_k (x).\quad \blacksquare
\end{eqnarray*}

\bigskip

Before performing the last step of our methodology, the usual
remark is in order. It is obvious that the above result is valid
only when $k$ is ``small". More precisely, the reader has to
remember that, if the starting alphabet has $m$ letters, then in
the above ECO matrix we must consider only the labels $(k_h )$
with $k\leq m$.

\bigskip

The last two theorems, stated as usual without proofs, conclude
our enumeration of $W_n ^{(m)}(1-11,1-12)$.

\begin{teor} The generating function for the number of words of $W_n
^{(m)}(1-11,1-12)$having $k$ distinct letters is
\begin{displaymath}
{m\choose k}(1+x)^k C_k (x).
\end{displaymath}
\end{teor}

\begin{teor}
The generating function of the sequence $|W_n ^{(m)}(1-11,1-12)|$
is
\begin{displaymath}
f_{\{1-11,1-12\}}^{(m)}=\sum_{k=0}^m {m\choose k}(1+x)^k C_k (x).
\end{displaymath}
\end{teor}

\section{Further work}

In our paper we have given complete results concerning the
enumeration of the classes of words $W(T)$, where $T$ consists of
two generalized patterns of length three of type $a-bc$ on the
alphabet $\{ 1,2\}$. Our main aim was to illustrate the soundness
of our methodology, as well as to use it to find some new counting
results. Of course, to complete the enumeration of $W(T)$ when $T$
is as above, many more cases should be considered. Following
essentially the same lines, it is possible to determine some
further generating functions, which we express below without
proof.

\begin{teor} The generating function of $|W_n ^{(m)}(1-11,1-21)|$
is
\begin{displaymath}
f_{\{1-11,1-21\}}^{(m)}=\sum_{k=0}^m {m\choose
k}x(1+x)^k\left(\prod_{i=0}^{k-2}\left((1+x)^{k-i}-1\right)\right)
\end{displaymath}
\end{teor}

\begin{teor} The generating function of $|W_n ^{(m)}(1-11,1-22)|$
is
\begin{displaymath}
f_{\{1-11,1-22\}}^{(m)}=\sum_{k=0}^m {m\choose
k}\frac{x^k(x+k)_k}{\prod_{i=1}^{k-1}(1-ix)}.
\end{displaymath}
\end{teor}

\begin{teor} The generating function of $|W_n ^{(m)}(2-11,1-22)|$
is
\begin{displaymath}
f_{\{2-11,1-22\}}^{(m)}=\sum_{k=0}^m \frac{(m)_k \cdot x^{k-1}
}{\left( \prod_{i=1}^{k-1} (1-ix)\right) \cdot (1-x)}.
\end{displaymath}
\end{teor}

\bigskip

The remaining cases seem to be more difficult to be dealt with; up
to Wilf-equivalence, they are the following:
\begin{eqnarray*}
W_n ^{(m)}(1-12,1-21)\qquad && \qquad W_n ^{(m)}(1-12,1-22)
\\ W_n ^{(m)}(1-12,2-11)\qquad && \qquad W_n ^{(m)}(1-12,2-12)
\\ W_n ^{(m)}(1-21,1-22)\qquad && \qquad W_n ^{(m)}(1-21,2-11).
\end{eqnarray*}

\end{document}